\font\largebf=cmbx10 scaled\magstep2
\font\tenmsa=msam10
\font\tenmsb=msbm10
\def\blackpen{\pdfliteral{.0 .0 .0 rg .0 .0 .0 RG}}
\def\redpen{\pdfliteral{.9 .0 .0 rg .9 .0 .0 RG}}
\long\def\red#1\endred{\redpen#1\blackpen}
\long\def\later#1{\red XXX\endred}

\def\all{\hbox{for all}}
\def\and{\hbox{and}}

\def\bigcupn{\bigcup\nolimits}
\def\bra#1#2{\langle#1,#2\rangle}
\def\Bra#1#2{\big\langle#1,#2\big\rangle}

\def\cite#1\endcite{[#1]}
\def\Cite#1\endcite{\big[#1\big]}

\def\dom{\hbox{\rm dom}}

\def\eps{\varepsilon}

\def\exs{\hbox{there exists}}

\def\f#1#2{{#1 \over #2}}

\def\half{\ts\f12}

\def\infn{\inf\nolimits}
\def\intr{\hbox{\rm int}}

\def\lr{\Longrightarrow}

\def\NQ{{\cal N}_q}
\def\on{\hbox{on}}
\def\PC{{\cal PC}}
\def\PCLSC{{\cal PCLSC}}

\def\PQ{{\cal P}_q}
\def\phi{\varphi}

\def\qed{\hfill\hbox{\tenmsa\char03}}
\def\qlr{\quad\lr\quad}

\def\quand{\quad\and\quad}
\def\r{\hbox{\rm I \hskip - .5em R}}
\def\r{\hbox{\tenmsb R}}
\def\rbar{\r\cup\{\infty\}}
\def\rbar{\,]{-}\infty,\infty]}

\def\rl{\Longleftarrow}

\def\rthalf{{\ts\f1{\sqrt2}}}

\long\def\slant#1\endslant{{\sl#1}}

\def\st{\hbox{such that}}
\def\supn{\sup\nolimits}

\def\toto{\ {\mathop{\hbox{\tenmsa\char19}}}\ }

\def\ts{\textstyle}

\def\xbra#1#2{\lfloor#1,#2\rfloor}
\def\Xbra#1#2{\big\lfloor#1,#2\big\rfloor}
\def\defSection#1{}
\def\defCorollary#1{}
\def\defDefinition#1{}
\def\defExample#1{}
\def\defLemma#1{}
\def\defNotation#1{}
\def\defProblem#1{}
\def\defRemark#1{}
\def\defTheorem#1{}
\def\locno#1{}
\def\meqno#1{\eqno(#1)}
\def\nmbr#1{}
\def\Proof{\medbreak\noindent{\bf Proof.}\enspace}
\def\Proo{{\bf Proof.}\enspace}
\def\Signoff{}
\def \INTsec{1}
\def \FENCHELsec{2}
\def \FENCHELthm{2.1}
\def \ABthm{2.2}
\def \EFnot{2.3}
\def \ABEPIthm{2.4}
\def \SSDsec{3}
\def \QPOSdef{3.1}
\def \QPOSone{1}
\def \QPOStwo{2}
\def \QPOSthree{3}
\def \QPOSfour{4}
\def \QPOSfive{5}
\def \Hex{3.2}
\def \EEex{3.3}
\def \NOTSSDex{3.4}
\def \POSFdef{3.5}
\def \POSFone{6}
\def \FATdef{3.6}
\def \FFATlem{3.7}
\def \PHlem{3.8}
\def \BCFNdef{3.9}
\def \BCFNone{7}
\def \TBCFNone{8}
\def \BCFNlem{3.10}
\def \PHIMAXthm{3.11}
\def \Jrem{3.12}
\def \TRANSlem{3.13}
\def \TBCBClem{3.14}
\def \ROOTlem{3.15}
\def \ROOTcor{3.16}
\def \Llem{3.17}
\def \FULLPHlem{3.18}
\def \SSDBsec{4}
\def \SSDBdef{4.1}
\def \SSDBthree{9}
\def \SSDBex{4.2}
\def \POSNEGthm{4.3}
\def \PTRANSthm{4.4}
\def \Grem{4.5}
\def \ABPOSNEGthm{4.6}
\def \ABPOSNEGcor{4.7}
\def \EEsec{5}
\def \EETRANSthm{5.1}
\def \REFLSUBDIFFrem{5.2}
\def \EPIthm{5.3}
\def \EPIone{10}
\def \EPItwo{11}
\def \COMBthm{5.4}
\def \EErem{5.5}
\def \EEone{12}
\def \MONsec{6}
\def \MAXthm{6.1}
\def \DDOMlem{6.2}
\def \SURthm{6.3}
\def \SURrem{6.4}
\def \SUMthm{6.5}
\def \SUMcor{6.6}
\def \RTRSUMcor{6.7}
\def \CQthm{6.8}
\def \HAMMsec{7}
\def \RHOTRANSthm{7.1}
\def \RHOSURthm{7.2}
\def \RHOSURcor{7.3}
\def \RRABthm{7.4}
\def \STINVthm{7.5}
\def \HAMMERthm{7.6}
\def \FIFTYrem{7.7}
\def \AB{1}
\def \FIFTY{2}
\def \BCW{3}
\def \BGW{4}
\def \BS{5}
\def \FITZ{6}
\def \MLSUR{7}
\def \MLT{8}
\def \PENOT{9}
\def \PZA{10}
\def \RTRFENCHEL{11}
\def \SUMS{12}
\def \MANDM{13}
\def \HBM{14}
\def \SSDMON{15}
\def \SZNZ{16}
\def \VZ{17}
\def \ZBOOK{18}
\def \ZFITZ{19}
\def \ZEIDLER{20}
\magnification 1200
\headline{\ifnum\folio=1
{\hfil{\largebf SSDB spaces and maximal monotonicity}\hfil}
\else\centerline{\rm {\bf SSDB spaces and maximal monotonicity}}\fi}
\medskip
\centerline{by}
\medskip
\centerline{Stephen Simons}
\bigskip
\centerline{\bf Abstract}
\medskip
\noindent
In this paper, we develop some of the theory of SSD spaces and SSDB spaces, and\break deduce some results on maximally monotone multifunctions on a reflexive Banach space.
\defSection \INTsec
\bigbreak
\centerline{\bf \INTsec\quad Introduction}
\medskip
\noindent In this paper, we develop enough of the theory of SSD spaces and SSDB spaces that we can obtain some significant results on maximally monotone multifunctions on a reflexive Banach space.   With a few minor additions, this is a written version of the lecture with the same title delivered at the IX ISORA meeting in Lima, Peru, in October 2009, and we will not attempt to give a comprehensive exposition of the theory of SSD spaces and SSDB spaces.   For this, we refer the reader to \cite\HBM\endcite\ and \cite\SSDMON\endcite, from which many of the proofs given here are taken.
\smallbreak
Many of the original results on maximally monotone multifunctions on a reflexive\break Banach space were obtained using Brouwer's fixed--point theorem either directly or indirectly.   The approach given here is based on convex analysis --- more specifically the Fenchel duality theorem.   In Section \FENCHELsec, we give the three versions of the Fenchel duality theorem that we will use in this paper.
\smallbreak
In Section \SSDsec, we introduce the concepts of \slant SSD space\endslant, \slant$q$--positive set\endslant\ and the functions $\Phi_\cdot$, which are the generalizations to SSD spaces of Fitzpatrick functions of monotone sets.   We also introduce the $q$--positive sets $\PQ(\cdot)$ determined by certain convex functions.   Despite the fact that this section uses the idea of conjugate function from convex analysis, the arguments used are essentially algebraic, apart from the disguised differentiability argument of Lemma \FFATlem.
\smallbreak
In Section \SSDBsec, we introduce the concept of \slant\ SSDB space\endslant, which is a SSD space with an appropriate Banach space structure.   The main results of this section are the two ``pos--neg theorems'', Theorem \POSNEGthm\ and \ABPOSNEGthm\ and the criterion for and properties of maximal $q$--positivity contained in Theorem \PTRANSthm.
\smallbreak
For the rest of this paper, we suppose that $E$ is a reflexive Banach space.   In Section \EEsec, we show how $E \times E^*$ can be considered as an SSDB space, and how Theorems \ABEPIthm\ and \PTRANSthm\ lead to results on maximal monotonicity.   Theorem \EETRANSthm\ and Theorem \COMBthm\ will be used in Section \MONsec. 
\smallbreak
Up to this point, our discussion of monotonicity has been in terms of monotone subsets of $E \times E^*$.   In Section \MONsec\ we move the emphasis to monotone multifunctions from $E$ into $E^*$.  In Theorem \SURthm, we show how Theorem \EETRANSthm\ gives Rockafellar's surjectivity theorem, and in Theorem \SUMthm, we show how Theorem \COMBthm\ gives sufficient conditions for the sum of maximally monotone multifunctions to be maximally monotone.
\smallbreak
In the final section of this paper, Section \HAMMsec, we show how the pos--neg theorem, Theorem \ABPOSNEGthm, together with the simple properties of the reflection maps $\rho_1,\rho_2\colon\ E \times E^* \to E \times E^*$ defined by $\rho_1(x,x^*) := (-x,x^*)$ and $\rho_2(x,x^*) := (x,-x^*)$ lead to an abstract Hammerstein theorem.
\smallbreak
The author would like to express his appreciation to Radu Ioan Bo\c t for some very constructive comments on the first version of this paper.  
\defSection\FENCHELsec
\bigbreak
\centerline{\bf \FENCHELsec\quad Versions of the Fenchel duality theorem}
\medskip
\noindent All vector spaces in the paper will be \slant real\endslant.   If $X$ is a nonzero vector space and $f\colon\ F \mapsto \rbar$ then we write  $\dom\,f := \{x \in F\colon\ f(x) \in \r\}$.   The set $\dom\,f$ is the \slant effective domain \endslant of $f$.   $f$ is said to be \slant proper\endslant\ if $\dom\,f \ne \emptyset$.   We write $\PC(X)$ for the set of proper convex functions from $X$ into $\rbar$.   If $X$ is a nonzero normed space we write $\PCLSC(X)$ for $\{f \in \PC(X)\colon\ f\ \hbox{is lower semicontinuous}\}$.   If $E$ and $F$ are nonzero vector spaces, $\bra\cdot\cdot\colon E \times F \to \r$ is a bilinear form and $f \in \PC(E)$ then the \slant Fenchel conjugate, $f^*$, of $f$ with respect to $\bra\cdot\cdot$\endslant\ is defined for $y \in F$ by
$$f^*(y) := \supn_{x \in E}\big[\bra{x}{y} - f(x)\big].$$
The commonest (but not the only) use of this notation is when $E$ is a normed space, $F$ is the dual, $E^*$, of $E$ and $\bra\cdot\cdot$ is the duality pairing. 
\medskip
We start off by stating a result that is an immediate consequence of
Rockafellar's version of the Fenchel duality theorem (see \cite\RTRFENCHEL, Theorem 1, p.\ 82--83\endcite\ for the original version and Z\u alinescu, \cite\ZBOOK, Theorem 2.8.7, p. 126--127\endcite\ for more general results):
\defTheorem \FENCHELthm
\medbreak
\noindent
{\bf Theorem \FENCHELthm.}\enspace\slant Let $F$ be a nonzero normed
space, $f\colon\ F \mapsto \rbar$ be proper and convex, $g\colon\
F \mapsto \r$ be convex and continuous, and $f + g \ge 0$ on $F$.
Then there exists $z^* \in F^*$ such that $f^*(z^*) + g^*(-z^*) \le 0$.
\endslant
\medskip
Theorem \ABthm\ below was first proved by Attouch--Brezis \big(this follows from \cite\AB, Corollary 2.3, pp.\ 131--132\endcite\big) -- there is a somewhat different proof in Simons, \cite\HBM, Theorem 15.1, p.\ 66\endcite, and a much more general result was established in Z\u alinescu, \cite\ZBOOK, Theorem 2.8.6, pp.\ 125--126\endcite.   We note that the result contained in Simons, \cite\HBM, Theorem 8.4, p.\ 46\endcite\ implies both Theorem \FENCHELthm\ and Theorem \ABthm, and that \cite\HBM, Theorem 7.4, p.\ 43\endcite\ gives a sharp lower bound on the possible values of $\|z^*\|$.
\defTheorem \ABthm
\medbreak
\noindent
{\bf Theorem \ABthm.}\enspace\slant Let $E$ be a nonzero Banach space, $f,\ g \in \PCLSC(E)$, $f + g \ge 0$ on $E$ and $\bigcup\nolimits_{\lambda > 0}\lambda\big[\dom\,f - \dom\,g\big]$ be a closed linear subspace of $E$.  Then there exists $z^* \in E^*$ such that $f^*(z^*) + g^*(-z^*) \le 0$.\endslant
\defNotation \EFnot
\medbreak
\noindent
{\bf Notation \EFnot.}\enspace If $E$ and $F$ are nonzero normed spaces, we define the dual of $E \times F$ to be $E^* \times F^*$ under the pairing $\Bra{(x,y)}{(x^*,y^*)} := \bra{x}{x^*} + \bra{y}{y^*}\quad\big((x,y) \in E \times F,\ (x^*,y^*) \in E^* \times F^*$\big).   We then define the \slant projection maps\endslant\ $\pi_1, \pi_2$ by $\pi_1(x,y) := x$ and $\pi_2(x,y) := y$.
\medbreak
We end this section with a bivariate generalization of Theorem \ABthm,   which was first proved in Simons--Z\u alinescu, \cite\SZNZ, Theorem 4.2, pp. 9--10\endcite.  There was a simpler proof given in Simons, \cite\HBM, Theorem 16.4, pp. 68--69\endcite.   The hypothesis of  Theorem \ABEPIthm\ is that $h(x,\cdot)$ is the ``inf--convolution'' of $f(x,\cdot)$ and $g(x,\cdot)$, and the conclusion is that $h^*(\cdot,y^*)$ is the ``exact inf--convolution'' of $f^*(\cdot,y^*)$ and $g^*(\cdot,y^*)$.
\defTheorem \ABEPIthm
\medbreak
\noindent
{\bf Theorem \ABEPIthm.}\enspace\slant Let $E$ and $F$ be nonzero Banach spaces, $f,\ g\colon\ E \times F \mapsto \rbar$ be proper, convex and lower semicontinuous, $\bigcup_{\lambda > 0}\lambda\big[\pi_1\,\dom\,f - \pi_1\,\dom\,g\big]$ be a closed linear subspace of $E$ and, for all $(x,y) \in E \times F$,
$$h(x,y) := \inf\big\{f(x,u) + g(x,v)\colon\ u,\,v \in F,\ u + v = y\big\} > -\infty.$$
Then, for all $(x^*,y^*) \in E^* \times F^* = (E \times F)^*$,
$$h^*(x^*,y^*) = \min\big\{f^*(s^*,y^*) + g^*(t^*,y^*)\colon\ s^*,\,t^* \in E^*,\ s^* + t^* = x^*\big\}.$$\endslant\par
\defSection \SSDsec
\medbreak
\centerline{\bf \SSDsec\quad SSD spaces}
\defDefinition \QPOSdef
\medbreak
\noindent
{\bf Definition \QPOSdef.}\enspace We will say that $\big(B,\xbra\cdot\cdot\big)$ is a \slant symmetrically self--dual space (SSD space)\endslant\  if $B$ is a nonzero real vector space and $\xbra\cdot\cdot\colon B \times B \to \r$ is a symmetric bilinear form.   In this case, we will \slant always\endslant\ write\quad $q(b) := \half\xbra{b}{b}$\quad($b \in B$).\quad   (``$q$'' stands for ``quadratic''.)
\smallskip
Now let $\big(B,\xbra\cdot\cdot\big)$ be an SSD space and $A \subset B$.   We say that $A$ is \slant$q$--positive\endslant\ if $A \ne \emptyset$ and
$$b,c \in A \lr q(b - c) \ge 0.$$
In this case, since $q(0) = 0$,
$$b \in A \lr \inf q(A - b) = 0.\meqno\QPOSone$$
We then define $\Phi_A\colon\ B \to \rbar$ by
$$\Phi_A(b) := \supn_A\big[\xbra\cdot{b} - q\big]\quad(b \in B).\meqno\QPOStwo$$
$\Phi_A$ is a generalization to SSD spaces of the ``Fiztpatrick function'' of a monotone set, which was originally introduced in \cite\FITZ\endcite\ in 1988, but lay dormant until it was rediscovered by  Mart\'\i nez-Legaz and Th\'era in \cite\MLT\endcite\ in 2001.
We note then that, for all $b \in B$,
$$\left.\eqalign{
\Phi_{A}(b) &= q(b) - \infn_{a \in A}\big[q(a) - \xbra{a}{b} + q(b)\big]\cr
&= q(b) - \infn_{a \in A}q(a - b) = q(b) - \inf q(A - b).}\right\}\meqno\QPOSthree$$
From (\QPOSone),
$$\Phi_A = q\ \on\ A.\meqno\QPOSfour$$
Thus $\Phi_A \in \PC(B)$.   We say that $A$ is \slant maximally $q$--positive\endslant\ if $A$ is $q$--positive and $A$ is not properly contained in any other $q$--positive set.   In this case, if $b \in B$ and $\inf q(A - b) \ge 0$ then clearly $b \in A$.   In other words,\qquad \big($b \in B \setminus A \lr \inf q(A - b) < 0$\big).\qquad From (\QPOSone),\qquad $\inf q(A - b) \le 0$\qquad and\qquad  \big($\inf q(A - b) = 0 \iff b \in A$\big).\qquad Thus, from (\QPOSthree)
$$\Phi_A \ge q\ \on\ B\quand\big(\Phi_A(b) = q(b) \iff b \in A\big).\meqno\QPOSfive$$
We make the elementary observation that if $b \in B$ and $q(b) \ge 0$ then the linear span $\r b$ of $\{b\}$ is $q$--positive.
\medbreak
We now give some examples of SSD spaces and their associated $q$--positive sets.   These examples are taken from \cite\HBM, pp.\ 79--80\endcite.   
\defExample \Hex 
\medbreak
\noindent
{\bf Example \Hex.}\enspace Let $B$ be a Hilbert space with inner product $(b,c) \mapsto \bra{b}{c}$.
\smallbreak
(a)\enspace If, for all $b,c \in B$, $\xbra{b}{c} := \bra{b}{c}$ then $q(b) = \half\|b\|^2$ and every nonempty subset of $B$ is $q$--positive.
\smallbreak
(b)\enspace If, for all $b,c \in B$, $\xbra{b}{c} := -\bra{b}{c}$ then $q(b) = -\half\|b\|^2$ and the $q$--positive sets are the singletons.
\smallbreak
(c)\enspace If $B = \r^3$ and 
$$\Xbra{(b_1,b_2,b_3)}{(c_1,c_2,c_3)} := b_1c_2 + b_2c_1 + b_3c_3,$$
then\quad$q(b_1,b_2,b_3) = b_1b_2 + \half b_3^2$.   Here, if $M$ is any nonempty monotone subset of $\r \times \r$ (in the obvious sense) then $M \times \r$ is a $q$--positive subset of $B$.   The set $\r(1,-1,2)$ is a $q$--positive subset of $B$ which is not contained in a set $M \times \r$ for any monotone subset of $\r \times \r$.   The helix $\big\{(\cos\theta,\sin\theta,\theta)\colon \theta \in \r\big\}$ is a $q$--positive subset of $B$, but if $0 < \lambda < 1$ then the helix $\big\{(\cos\theta,\sin\theta,\lambda\theta)\colon \theta \in \r\big\}$ is not.
\defExample \EEex 
\medbreak
\noindent
{\bf Example \EEex.}\enspace Let $E$ be a nonzero Banach space and $B := E \times E^*$.   For all $(x,x^*)$ and $(y,y^*) \in B$, we set\quad $\Xbra{(x,x^*)}{(y,y^*)} := \bra{x}{y^*} + \bra{y}{x^*}$.\quad  Then $\big(B,\xbra\cdot\cdot\big)$ is an SSD space with\quad $q(x,x^*) = \half\big[\bra{x}{x^*} + \bra{x}{x^*}\big] = \bra{x}{x^*}$.\quad   Consequently, if\quad $(x,x^*), (y,y^*) \in B$\quad then\quad$\bra{x - y}{x^* - y^*} = q(x - y,x^* - y^*) = q\big((x,x^*) - (y,y^*)\big)$.\quad   
Thus if $A \subset B$ then $A$ is $q$--positive exactly when $A$ is a nonempty monotone subset of $B$ in the usual sense, and $A$ is maximally $q$--positive exactly when $A$ is a maximally monotone subset of $B$ in the usual sense.   We point out that any finite dimensional SSD space of the form described here must have \slant even\endslant\ dimension.  Thus cases of Example \Hex\ with finite odd dimension cannot be of this form.
\defExample \NOTSSDex
\medbreak
\noindent
{\bf Example \NOTSSDex.}\enspace $\big(\r^3,\xbra\cdot\cdot\big)$ is \slant not\endslant\ an SSD space with
$$\Xbra{(b_1,b_2,b_3)}{(c_1,c_2,c_3)} := b_1c_2 + b_2c_3 + b_3c_1.$$
(The bilinear form $\xbra\cdot\cdot$ is not symmetric.)
\defDefinition \POSFdef
\medbreak
\noindent
{\bf Definition \POSFdef.}\enspace Let $\big(B,\xbra\cdot\cdot\big)$ be an SSD space.   If $f \in \PC(B)$ and $f \ge q$ on $B$, we write
$$\PQ(f) := \big\{b \in B\colon\ f(b) = q(b)\big\}.$$
We then note from  (\QPOSfive) that
$$\hbox{\sl if}\ A\ \hbox{\sl is maximally}\ q\hbox{\sl--positive then}\ A = \PQ(\Phi_A).\meqno\POSFone$$
If $g \in \PC(B)$ and  $g \ge -q$ on $B$ then we write $\NQ(g) := \big\{b \in B\colon\ g(b) = -q(b)\big\}$.
\medbreak
We now introduce the concept of \slant intrinsic conjugate\endslant\ for an SSD space.
\defDefinition \FATdef
\medbreak
\noindent
{\bf Definition \FATdef.}\enspace If $\big(B,\xbra\cdot\cdot\big)$ is an SSD space and $f \in \PC(B)$, we write $f^@$ for the Fenchel conjugate of $f$ with respect to the pairing $\xbra\cdot\cdot$, that is to say,
$$\all\ c \in B,\qquad f^@(c) := \supn_B\big[\xbra{\cdot}{c} - f\big].$$
\par
The concepts introduced in Definitions \POSFdef\ and \FATdef\ are related in our next result, which uses a disguised differentiability argument.   
\defLemma \FFATlem
\medbreak
\noindent
{\bf Lemma \FFATlem.}\enspace\slant Let $\big(B,\xbra\cdot\cdot\big)$ be an SSD space,\quad $f \in \PC(B)$\quand  $f \ge q$ on $B$.\quad   Then $f^@ = q\ \on\ \PQ(f)$.
\endslant
\Proof Let $a \in \PQ(f)$.   Let $\lambda \in \,]0,1[\,$.   For simplicity in writing, let $\mu := 1 - \lambda \in \,]0,1[\,$.   Then, for all $b \in B$,
$$\eqalign{\lambda^2q(b) + \lambda\mu\xbra{b}{a} + \mu^2q(a) &= q\big(\lambda b + \mu a\big) \le f(\lambda b + \mu a)\cr
&\le \lambda f(b) + \mu f(a) = \lambda f(b) + \mu q(a).}$$
Thus\quad $\lambda^2q(b) + \lambda\mu\xbra{b}{a} \le \lambda f(b) + \lambda\mu q(a)$.\quad Dividing by $\lambda$ and letting $\lambda \to 0$, we have\quad $\xbra{b}{a} \le f(b) + q(a)$,\quad that is to say\quad$\xbra{a}{b} - f(b) \le q(a)$,\quad and, taking the supremum over $b$, $f^@(a) \le q(a)$.   On the other hand,\quad $f^@(a) \ge \xbra{a}{a} - f(a) = 2q(a) - q(a) = q(a)$,\quad completing the proof of Lemma \FFATlem.\qed
\medbreak
The next result gives a basic property of ${\Phi_A}^@$.
\defLemma \PHlem
\medbreak
\noindent
{\bf Lemma \PHlem.}\enspace\slant Let $\big(B,\xbra\cdot\cdot\big)$ be an SSD space and $A$ be a nonempty $q$--positive subset of $B$.   Then\quad ${\Phi_A}^@ \ge  \Phi_A$ on $B$.\endslant
\Proof Let $c \in B$.   Then, from (\QPOSfour), 
$${\Phi_A}^@(c) = \supn_B\big[\xbra{c}\cdot - \Phi_A\big]
\ge \supn_A\big[\xbra{c}\cdot - \Phi_A\big]
= \supn_A\big[\xbra{c}\cdot - q\big] = \Phi_A(c).\eqno\qed$$\par
We now introduce the important concepts of ``BC--function'' and ``TBC--function''.
\defDefinition \BCFNdef
\medbreak
\noindent
{\bf Definition \BCFNdef.}\enspace Let $\big(B,\xbra\cdot\cdot\big)$ be a SSD space and $f,g \in \PC(B)$.   We say that $f$ is a\slant\ BC--function\endslant\ if
$$b \in B \qlr f^@(b) \ge f(b) \ge q(b).\meqno\BCFNone$$
``BC'' stands for ``bigger conjugate''.  We say that $g$ is a\slant\ TBC--function\endslant\ if
$$b \in B \qlr g^@(-b) \ge g(b) \ge -q(b).\meqno\TBCFNone$$
``TBC'' stands for ``twisted bigger conjugate''.   Of course, $g$ is a TBC--function if, and only if, $g$ is a BC--function with respect to the SSD space $\big(B,-\xbra\cdot\cdot\big)$ 
\defLemma \BCFNlem
\medbreak
\noindent
{\bf Lemma \BCFNlem.}\enspace\slant Let $\big(B,\xbra\cdot\cdot\big)$ be a SSD space and $f \in \PC(B)$ be a BC--function.   Then 
$$\PQ\big(f^@\big) = \PQ(f).$$\endslant
\Proo This follows from Lemma \FFATlem\ and the inclusion\quad $\PQ\big(f^@\big) \subset \PQ(f)$,\quad which is\break immediate from (\BCFNone).\qed
\defTheorem  \PHIMAXthm
\medbreak
\noindent
{\bf Theorem  \PHIMAXthm.}\enspace\slant Let $\big(B,\xbra\cdot\cdot\big)$ be an SSD space and $A$ be a maximally $q$--positive subset of $B$.   Then\quad ${\Phi_A}$ is a BC--function\quand $\PQ\big({\Phi_A}^@\big) = \PQ\big(\Phi_A\big) = A$.\endslant
\Proof The first assertion follows from Lemma \PHlem\ and (\QPOSfive), and the second assertion is immediate from Lemma \BCFNlem\ and (\POSFone).\qed   
\defRemark \Jrem
\medbreak
\noindent
{\bf Remark \Jrem.}\enspace We shall see by combining (\SSDBthree) and Remark \EErem\ that there may exist a function that is both a BC--function and a TBC--function but not of the form $\Phi_A$ for any nonempty $q$--positive set $A$. \medbreak
We now give two computational results.   In the first of these, we investigate the ``translation'' of a BC--function.   This result will be used in the ``pos--neg'' theorems, Theorem \POSNEGthm\ and Theorem \ABPOSNEGthm. 
\defLemma \TRANSlem
\medbreak
\noindent
{\bf Lemma \TRANSlem.}\enspace\slant Let $\big(B,\xbra\cdot\cdot\big)$ be a SSD space, $f \in \PC(B)$ be a BC--function and $c \in B$.   We define\quad $f_c\in \PC(B)$\quad by\quad $f_c := f(\cdot + c) - \xbra{\cdot}{c} - q(c)$.\quad  Then $f_c$ is a BC--function,\quad $\dom\,f_c = \dom\,f - c$\quad and\quad $\PQ(f_c) = \PQ(f) - c$.\endslant
\Proof For all $b \in B$,
$$\leqalignno{{f_c}^@(b)
&= \supn_{d \in B}\big[\xbra{d}{b} + \xbra{d}{c} + q(c) - f(d + c)\big]\cr
&= \supn_{e \in B}\big[\xbra{e - c}{b + c} + q(c) - f(e)\big]\cr
&= \supn_{e \in B}\big[\xbra{e}{b + c} - \xbra{c}{b} - f(e)\big] - q(c)\cr
&= f^@(b + c) - \xbra{c}{b} - q(c).}$$
It follows from (\BCFNone) that ${f_c}^@(b) \ge f(b + c) - \xbra{b}{c} - q(c) = f_c(b)$ and
$$f_c(b) = f(b + c) - \xbra{b}{c} - q(c) \ge q(b + c) - \xbra{b}{c} - q(c) = q(b).$$
Consequently, $f_c$ is a BC--function.   It is obvious that $\dom\,f_c = \dom\,f - c$.   Further, since
$$\eqalign{b \in \PQ(f_c) &\iff f(b + c) - \xbra{b}{c} - q(c) = q(b)\cr
&\iff f(b + c) = q(b + c) \iff b + c \in \PQ(f),}$$
we have $\PQ(f_c) = \PQ(f) - c$, as required.\qed
\medbreak
In the second computational result, we show how we can always obtain a TBC--function from a BC--function by an appropriate linear transformation.
\defLemma \TBCBClem
\medbreak
\noindent
{\bf Lemma \TBCBClem.}\enspace\slant Let $\big(B,\xbra\cdot\cdot\big)$ be an SSD space, $\rho\colon B \to B$ be a linear surjection such that,\quad for all $b,c \in B$, $\Xbra{\rho(b)}{\rho(c)} = \xbra{b}{-c}$, \quad and $g \in \PC(B)$ be a BC--function.  Then\quad $g \circ \rho$ and $g \circ (-\rho)$ are TBC-- functions,\quad $\rho\,\dom\,(g \circ \rho) = \dom\,g$,\quad $\rho\NQ(g \circ \rho) = \PQ(g)$, $-\rho\,\dom\,(g \circ (-\rho)) = \dom\,g$ \quand $-\rho\NQ(g \circ (-\rho)) = \PQ(g)$.\endslant
\Proof We give the proof for $g \circ \rho$.   The proof for $g \circ (-\rho)$ is similar.   For all $c \in B$,
$$\eqalign{{(g \circ \rho)}^@(-c) &= \supn_{b \in B}\big[\xbra{b}{-c} - g\big(\rho(b)\big)\big] = \supn_{b \in B}\big[\xbra{\rho(b)}{\rho(c)} - g\big(\rho(b)\big)\big]\cr
&= \supn_{d \in B}\big[\xbra{d}{\rho(c)} - g(d)\big] = g^@\big(\rho(c)\big) \ge g\big(\rho(c)\big)\cr 
&\ge q\big(\rho(c)\big) = \half\xbra{\rho(c)}{\rho(c)} = \half\xbra{c}{-c} = -q(c).}$$
Thus $g \circ \rho$ is a TBC-- function.   Further,
$$\eqalign{b \in \rho\,\dom\,(g \circ \rho)
&\iff \exs\ c \in B\ \st\ g\big(\rho(c)\big) \in \r\ \and\ b = \rho(c)\cr
&\iff \exs\ d \in B\ \st\ g(d) \in \r\ \and\ b = d\cr
&\iff g(b) \in \r \iff b \in \dom\,g,}$$
from which\quad $\rho\,\dom\,(g \circ \rho) = \dom\,g$,\quad and
$$\eqalign{b \in \rho\,\NQ(g \circ \rho)
&\iff\ \exs\ c \in B\ \st\ g\big(\rho(c)\big) = -q(c)\ \and\ b = \rho(c)\cr
&\iff\ \exs\ c \in B\ \st\ g\big(\rho(c)\big) = q\big(\rho(c)\big)\ \and\ b = \rho(c)\cr
&\iff\ \exs\ d \in B\ \st\ g(d) = q(d)\ \and\ b = d\cr
&\iff b \in \PQ(g),}$$
from which\quad $\rho\NQ(g \circ \rho) = \PQ(g)$.\qed
\medbreak
Lemma \ROOTlem\ is a subtler property of the bilinear form $q$.   In the situation of Example \EEex, Corollary \ROOTcor\ gives us Voisei--Z\u alinescu \cite\VZ, Proposition 1\endcite.    
\defLemma \ROOTlem
\medbreak
\noindent
{\bf Lemma \ROOTlem.}\enspace\slant Let $\big(B,\xbra\cdot\cdot\big)$ be an SSD space, $f \in \PC(B)$, $f \ge q$ on $B$ and $b,c \in B$.   Then
$$-q(b - c) \le \Big[\sqrt{(f - q)(b)} + \sqrt{(f - q)(c)}\Big]^2.$$\endslant
\Proof See \cite\SSDMON, Lemma 2.6\endcite.\qed 
\defCorollary \ROOTcor
\medbreak
\noindent
{\bf Corollary \ROOTcor.}\enspace\slant Let $\big(B,\xbra\cdot\cdot\big)$ be an SSD space, $f \in \PC(B)$, $f \ge q$ on $B$ and $b,c \in B$.   Then
$$-q(b - c) \le 2(f - q)(b) + 2(f - q)(c).$$\endslant
\Proo This is immediate from Lemma \ROOTlem\ and the Cauchy--Schwarz inequality.\qed 
\medbreak
Lemma \Llem\ is suggested by Burachik--Svaiter, \cite\BS, Theorem 3.1, pp. 2381--2382\endcite\ and Penot, \cite\PENOT, Proposition 4(h)$\lr$(a), pp. 860--861\endcite.   Lemma \Llem\ will be used in Theorem \PTRANSthm.
\defLemma \Llem
\medbreak
\noindent
{\bf Lemma \Llem.}\enspace\slant Let $\big(B,\xbra\cdot\cdot\big)$ be an SSD space, $f \in \PC(B)$, $f \ge q$ on $B$ and $\PQ(f) \ne \emptyset$.   Then $\PQ(f)$ is a $q$--positive subset of $B$.\endslant
\Proof This is immediate from Lemma \ROOTlem, or Corollary \ROOTcor, or by using the fact that $q$ satisfies the \slant parallelogram law\endslant:\quad $x,y \in B \lr q(x) + q(y) = \half q(x + y) + \half q(x - y)$.\qed
\medbreak
The final result of this section can be thought of as a completion of Lemma \PHlem.   It will not be used in the sequel.
\defLemma \FULLPHlem
\medbreak
\noindent
{\bf Lemma \FULLPHlem.}\enspace\slant Let $\big(B,\xbra\cdot\cdot\big)$ be an SSD space and $A$ be a nonempty $q$--positive subset of $B$.   Then\quad ${\Phi_A}^@ \le q$ on $A$, \quad  ${\Phi_A}^@ \ge  q$ on $B$,\quand ${\Phi_A}^{@@} = \Phi_A$ on $B$.\endslant
\Proof See \cite\SSDMON, Lemma 2.11\endcite.\qed
\defSection \SSDBsec
\medbreak
\centerline{\bf \SSDBsec\quad SSDB spaces}
\medskip
\noindent
We now introduce the SSDB spaces, a subclass of the class of SSD spaces.   Our treatment of the SSD spaces in Section \SSDsec\ has been essentially nontopological.   The additional norm structure of the SSDB spaces is essentially what makes maximally monotone multifunctions on a reflexive Banach space much more tractable than those on a general Banach space.   We will return to this issue when we consider Example \EEex\ in Section \EEsec.
\defDefinition \SSDBdef
\medbreak
\noindent
{\bf Definition \SSDBdef.}\enspace We will say that $\big(B,\xbra\cdot\cdot,\|\cdot\|\big)$ is a \slant symmetrically self--dual Banach space (SSDB space)\endslant\ if $\big(B,\xbra\cdot\cdot\big)$ is a SSD space, $\big(B,\|\cdot\|\big)$ is a Banach space, and there exists a linear isometry $\iota$ from $B$ onto $B^*$ such that, for all $b,c \in B$, $\Bra{b}{\iota(c)} = \xbra{b}{c}$.   We note then that, for all $f \in \PC(B)$ and $c \in B$,\quad $f^@(c) := \supn_B\big[\xbra\cdot{c} - f\big] = \supn_B\big[\bra\cdot{\iota(c)} - f\big] = f^*\big(\iota(c)\big)$,\quad that is to say\quad $f^@ = f^* \circ \iota$.   It is easy to see that the quadratic form $q$ is continuous and, for all $b \in B$,\quad $|q(b)|  = \half\big|\xbra bb\big| \le  \half\|b\|^2$.
\smallbreak
Let $g_0 := \half\|\cdot\|^2$ on $B$.   Then, for all $b \in B$,
$${g_0}^@(b) = {g_0}^*\big(\iota(b)\big) = \half\|\iota(b)\|^2 = \half\|b\|^2 = g_0(b).$$
It follows that
$$g_0\ \hbox{is both a BC--function and a TBC--function.}\meqno\SSDBthree$$\par
%
\defExample \SSDBex 
\medbreak
\noindent
{\bf Examples \SSDBex.}\enspace
(a)\enspace In Example \Hex(a), $\big(B,\xbra\cdot\cdot,\|\cdot\|\big)$ is a SSDB space under the Hilbert space norm, $q = g_0$ and $\NQ(g_0) = \{0\}$.   \big(We recall that the set $\NQ(g_0)$ was defined in Definition \POSFdef.\big)
\smallbreak
(b)\enspace  In Example \Hex(b), $\big(B,\xbra\cdot\cdot,\|\cdot\|\big)$ is a SSDB space under the Hilbert space norm, $q = -g_0$ and $\NQ(g_0) = B$.
\smallbreak
(c)\enspace In Example \Hex(c), $\big(B,\xbra\cdot\cdot,\|\cdot\|\big)$ is a SSDB space under the Euclidean norm and
$$\NQ(g_0) = \{(b_1,b_2,b_3) \in B\colon\ b_1 + b_2 = 0,\ b_3 = 0\}.$$
\smallbreak
Our next result is the ``pos--neg theorem of Rockafellar type''.   This will be used indirectly in Theorem \EETRANSthm, which will be used in turn indirectly in our proof of Rockafellar's surjectivity theorem, Theorem \SURthm, and the sum theorem, Theorem \SUMthm.   This result first appeared in Simons, \cite\HBM, Theorem 19.16, p.\ 83\endcite. 
\defTheorem \POSNEGthm
\medbreak
\noindent
{\bf Theorem \POSNEGthm.}\enspace\slant Let $\big(B,\xbra\cdot\cdot,\|\cdot\|\big)$ be a SSDB space, $f \in \PC(B)$ be a BC--function and $g\colon B \to \r$ be a  continuous TBC--function.  Then\quad$\PQ(f) - \NQ(g) = B$.\endslant
\Proof Let $c$ be an arbitrary element of $B$.   From Lemma \TRANSlem,   $f_c$ is a BC--function and so, using (\BCFNone) and (\TBCFNone),
$$b \in B \qlr f_c(b) + g(b) \ge q(b) - q(b) = 0.$$
Rockafellar's version of the Fenchel duality theorem, Theorem \FENCHELthm, and the surjectivity of $\iota$ now give $b \in B$ such that\quad ${f_c}^*\big(\iota(b)\big) + g^*\big({-}\iota(b)) \le 0$,\quad that is to say\quad ${f_c}^@(b) + g^@(-b) \le 0$.\quad  From (\BCFNone) and (\TBCFNone) again,\quad $f_c(b) + g(b) \le 0 = q(b) - q(b)$.\quad   From (\BCFNone) and (\TBCFNone) for a third time,\quad $f_c(b) = q(b)$ and $g(b) = -q(b)$,\quad that is to say, using Lemma \TRANSlem,\quad $b \in \PQ(f_c) = \PQ(f) - c$ and also $b \in \NQ(g)$.\quad   But then\quad $c = (c + b) - b \in \PQ(f) - \NQ(g)$.\quad   This completes the proof of Theorem \POSNEGthm.\qed
\defTheorem \PTRANSthm
\medbreak
\noindent
{\bf Theorem \PTRANSthm.}\enspace\slant Let $\big(B,\xbra\cdot\cdot,\|\cdot\|\big)$ be a SSDB space.   Then:
\smallbreak
\noindent
{\rm(a)}\enspace Let $f \in \PC(B)$ be a BC--function.   Then\quad  $\PQ(f) - \NQ(g_0) = B$.
\smallbreak
\noindent
{\rm(b)}\enspace Let $A$ be a $q$--positive subset of $B$ and\quad $A - \NQ(g_0) = B$.\quad Then $A$ is maximally $q$--positive.
\smallbreak\noindent
{\rm(c)}\enspace Let $f \in \PC(B)$ be a BC--function.   Then\quad $\PQ(f)$\quad is maximally $q$--positive.   Furthermore,\quad $\PQ(f^@) = \PQ(f)$.
\smallbreak
\noindent
{\rm(d)}\enspace Let $A$ be a $q$--positive subset of $B$.   Then\quad $A$ is maximally $q$--positive if, and only if\quad $A - \NQ(g_0) = B$. 
\endslant
\Proof(a) This follows from (\SSDBthree) and Theorem \POSNEGthm.
\smallbreak
(b) Suppose that $b \in B$ and\quad $A \cup \{b\}$\quad is $q$--positive.   By hypothesis, there exists $a \in A$ such that\quad $a - b \in \NQ(g_0)$.\quad   Thus\quad $\half\|a - b\|^2 = -q(a - b)$.\quad   Since\quad $A \cup \{b\}$\quad is $q$--positive,\quad $q(a - b) \ge 0$,\quad and so\quad $\half\|a - b\|^2 \le 0$.\quad   Thus\quad $b = a \in A$.
\smallbreak
(c)\enspace From (a) and Lemma \Llem, $\PQ(f)$ is nonempty and $q$--positive.   (c) now follows from (a), (b) and Lemma \BCFNlem.
\smallbreak
(d)\enspace We have already established $(\rl)$ in (b).   Suppose, conversely, that $A$ is maximally $q$--positive.   It follows from Theorem \PHIMAXthm\ that $\Phi_A$ is a BC--function and $\PQ(\Phi_A) = A$.   Thus, from (a), $A - \NQ(g_0) = \PQ(\Phi_A) - \NQ(g_0) = B$.\qed 
\defRemark \Grem
\medbreak
\noindent
{\bf Remark \Grem.}\enspace We note from (\SSDBthree) and Theorem \PTRANSthm(c) that $\PQ(g_0)$ is maximally $q$--positive and $\PQ({g_0}^@) = \PQ(g_0)$.   Now suppose that $a \in \PQ(g_0)$ and $b \in B$.   From the Fenchel--Young inequality,\quad $\xbra{a}{b} - q(a) = \xbra{a}{b} - {g_0}^@(a) \le g_0(b)$.\quad Taking the supremum over $a \in  \PQ(g_0)$, we deduce that $\Phi_{\PQ(g_0)} \le g_0$ on $B$.   In particular, $\dom\,\Phi_{\PQ(g_0)} = B$. 
\medbreak
The next result in this section is the ``pos-neg theorem of Attouch--Brezis type''.   A stronger result was proved in Simons, \cite\HBM, Theorem 21.12, p.\ 93\endcite.
\defTheorem \ABPOSNEGthm
\medbreak
\noindent
{\bf Theorem \ABPOSNEGthm.}\enspace\slant Let $\big(B,\xbra\cdot\cdot,\|\cdot\|\big)$ be a SSDB space, $f \in \PCLSC(B)$ be a BC--function, and $g \in \PCLSC(B)$ be a TBC--function.   Then\quad $\dom\,f - \dom\,g = B \iff\PQ(f) - \NQ(g) = B$.\endslant
\Proof Since the implication ($\rl$) is trivial, we only have to prove ($\lr$).   Let $c$ be an arbitrary element of $B$.   Then, from Lemma \TRANSlem,
$$\dom\,f_c - \dom\,g = \dom\,f - c - \dom\,g = B - c = B,$$
and so\quad $\bigcup\nolimits_{\lambda > 0}\lambda\big[\dom\,f_c - \dom\,g\big] = B$.\quad   From Lemma \TRANSlem,   $f_c$ is a BC--function and so, using (\BCFNone) and (\TBCFNone),
$$b \in B \qlr f_c(b) + g(b) \ge q(b) - q(b) = 0.$$
The Attouch--Brezis theorem, Theorem \ABthm,  and the surjectivity of $\iota$ now give $b \in B$ such that\quad ${f_c}^*\big(\iota(b)\big) + g^*\big({-}\iota(b)) \le 0$,\quad that is to say\quad ${f_c}^@(b) + g^@(-b) \le 0$.\quad From (\BCFNone) and (\TBCFNone) again,\quad $f_c(b) + g(b) \le 0 = q(b) - q(b)$.\quad   From (\BCFNone) and (\TBCFNone) for a third time,\quad  $f_c(b) = q(b)$ and $g(b) = -q(b)$,\quad that is to say, using Lemma \TRANSlem,\quad $b \in \PQ(f_c) = \PQ(f) - c$ and also $b \in \NQ(g)$.\quad   But then\quad $c = (c + b) - b \in \PQ(f) - \NQ(g)$.\quad  This completes the proof of Theorem \ABPOSNEGthm.\qed
\medbreak
Corollary \ABPOSNEGcor\ below is the form in which we will actually apply Theorem \ABPOSNEGthm\ to abstract Hammerstein theorems in Theorem \RHOTRANSthm. 
\defCorollary \ABPOSNEGcor
\medbreak
\noindent
{\bf Corollary \ABPOSNEGcor.}\enspace\slant Let $\big(B,\xbra\cdot\cdot,\|\cdot\|\big)$ be a SSDB space, $f,g \in \PCLSC(B)$ and $f,g$ be BC--functions.   Suppose, further, that $\rho\colon B \to B$ is a continuous linear bijection such that,\quad for all $b,c \in B$, $\xbra{\rho(b)}{\rho(c)} = \xbra{b}{-c}$.\quad   Then
$$\dom\,f - \rho^{-1}\dom\,g  = B \iff \PQ(f)  - \rho^{-1}\PQ(g) = B$$
and
$$\dom\,f + \rho^{-1}\dom\,g  = B \iff \PQ(f)  + \rho^{-1}\PQ(g) = B.$$\endslant
\Proof This is immediate from Theorem \ABPOSNEGthm\ and Lemma \TBCBClem.\qed
\defSection \EEsec
\bigbreak
\centerline{\bf \EEsec\quad BC--functions on $E \times E^*$}
\medskip
\noindent
From now on, $E$ is a nonzero reflexive Banach space and $E^*$ is its topological dual space.  As observed in Example \EEex, $\big(E \times E^*,\xbra\cdot\cdot\big)$ is an SSD space with $q(x,x^*) = \bra{x}{x^*}$, and if $A \subset E \times E^*$ then $A$ is $q$--positive exactly when $A$ is a nonempty monotone set in the usual sense, and $A$ is maximally $q$--positive exactly when $A$ is a maximally monotone set of in the usual sense. We norm $E \times E^*$ by $\big\|(x,x^*)\big\| := \sqrt{\|x\|^2 + \|x^*\|^2}$.   Then
$$\big(E \times E^*,\|\cdot\|\big)^* = (E^* \times E,\|\cdot\|\big),$$
under the duality $\Bra{(x,x^*)}{(y^*,y)} := \bra{x}{y^*} + \bra{y}{x^*}$.   Combining this with the formula in Example \EEex, we have\quad 
$\Bra{(x,x^*)}{(y^*,y)} = \Xbra{(x,x^*)}{(y,y^*)}$,\quad and so\quad $\iota(y,y^*) = (y^*,y)$.   It is easily seen from these relationships that $\big(E \times E^*,\xbra\cdot\cdot,\|\cdot\|\big)$ is an SSDB space.   Further,\quad$(x,x^*) \in \PQ(g_0) \iff \half\|x\|^2 + \half\|x^*\|^2 = \bra{x}{x^*} \iff (x,x^*) \in G(J),$\quad and\quad$(x,x^*) \in \NQ(g_0) \iff \half\|x\|^2 + \half\|x^*\|^2 = -\bra{x}{x^*} \iff (x,x^*) \in G(-J),$\quad where\quad $J\colon E \toto E^*$\quad is the duality map.   Thus we obtain the following fundamental\break result:
\defTheorem \EETRANSthm
\medbreak
\noindent
{\bf Theorem \EETRANSthm.}\enspace\slant Let $E$ be a nonzero reflexive Banach space.   Then:\smallbreak
\noindent
{\rm(a)}\enspace   If $f \in \PC(E \times E^*)$ and $f$ is a BC--function then\quad $\PQ(f)$ is maximally monotone,\quad $\PQ\big(f^@\big) = \PQ(f)$\quand $\PQ(f) - G(-J) = B$.
\smallbreak
\noindent
{\rm(b)}\enspace If $A$ is a monotone subset of $E \times E^*$ then $A$ is maximally monotone if, and only if, $A - G(-J) = E \times E^*$. 
\endslant
\Proof (a) is immediate from Theorem \PTRANSthm(a,c), and (b) is immediate from Theorem \PTRANSthm(d).\qed 
\defRemark \REFLSUBDIFFrem
\medbreak
\noindent
{\bf Remark \REFLSUBDIFFrem.}\enspace If $f \in \PCLSC(E)$ and, for all $(x,x^*) \in B$, $h(x,x^*) := f(x) + f^*(x^*)$ then the Fenchel--Moreau theorem and the Fenchel--Young inequality imply that $h$ is a BC--function on $E \times E^*$.   Since $\PQ(h) = G(\partial f)$, Theorem \EETRANSthm(a) gives us that $\partial f$ is maximally monotone.  (Remember that we are assuming that $E$ is reflexive!)  
\medbreak
We note that the ``partial episum theorem'' of Theorem \EPIthm\ below can also be deduced from Penot--Z\u alinescu, \cite\PZA,  Corollary 3.7\endcite.
\defTheorem \EPIthm
\medbreak
\noindent
{\bf Theorem \EPIthm.}\enspace\slant Let $E$ be a nonzero reflexive Banach space, $f,\ g \in\ \PCLSC(E \times E^*)$ be BC--functions, $\bigcup_{\lambda > 0}\lambda\big[\pi_1\,\dom\,f - \pi_1\,\dom\,g\big]$ be a closed linear subspace of $E$ and, for all $(x,x^*) \in E \times E^*$,
\locno\EPIone
$$h(x,x^*) := \inf\big\{f(x,s^*) + g(x,t^*)\colon\ s^*,\,t^* \in E^*,\ s^* + t^* = x^*\big\}.\eqno(\EPIone)$$
Then $h$ is a BC--function,\quad $\PQ(h^@) = \big\{(x,s^* + t^*)\colon\ (x,s^*) \in \PQ(f^@),\ (x,t^*) \in \PQ(g^@)\big\}$,\quad and
$$\PQ(h) = \big\{(x,s^* + t^*)\colon\ (x,s^*) \in \PQ(f),\ (x,t^*) \in \PQ(g)\big\}.$$
\endslant
\Proof Since $f \ge q$ and $g \ge q$ on $E \times E^*$, (\EPIone) implies that, for all $(x,x^*) \in E \times E^*$,
$$h(x,x^*) \ge \inf\big\{\bra{x}{s^*} + \bra{x}{t^*}\colon\ s^*,\,t^* \in E^*,\ s^* + t^* = x^*\big\} = \bra{x}{x^*} = q(x,x^*),$$
and then Theorem \ABEPIthm\ and the fact that $f^@ \ge f$ and $g^@ \ge g$ on $E \times E^*$ give
\locno\EPItwo
$$\eqalignno{h^@(x,x^*) = h^*(x^*,x)
&= \min\big\{f^*(s^*,x) + g^*(t^*,x)\colon\ s^*,\,t^* \in E^*,\ s^* + t^* = x^*\big\}\cr
&= \min\big\{f^@(x,s^*) + g^@(x,t^*)\colon\ s^*,\,t^* \in E^*,\ s^* + t^* = x^*\big\}&(\EPItwo)\cr
&\ge \inf\big\{f(x,s^*) + g(x,t^*)\colon\ s^*,\,t^* \in E^*,\ s^* + t^* = x^*\big\} = h(x,x^*).}$$
Thus $h$ is a BC--function, the required characterization of $\PQ(h^@)$ is immediate from (\EPItwo), and the final assertion now follows from three applications of Theorem \EETRANSthm(a).\qed
\defTheorem \COMBthm
\medbreak
\noindent
{\bf Theorem \COMBthm.}\enspace\slant Let $E$ be a nonzero reflexive Banach space, $f,\ g \in\ \PCLSC(E \times E^*)$ be BC--functions and $\bigcup_{\lambda > 0}\lambda\big[\pi_1\,\dom\,f - \pi_1\,\dom\,g\big]$ be a closed linear subspace of $E$.   Then
$$\big\{(x,s^* + t^*)\colon\ (x,s^*) \in \PQ(f),\ (x,t^*) \in \PQ(g)\big\}$$
is a maximally monotone subset of $E \times E^*$.\endslant
\Proof This is immediate from Theorem \EETRANSthm(a) and Theorem \EPIthm.\qed
\defRemark \EErem
\medbreak
\noindent
{\bf Remark \EErem.}\enspace The question arises naturally as to whether there can exist a $q$--positive subset $A$ of $B$ such that $\Phi_A = g_0$.   It would then follow from (\QPOSfour) that $g_0 = q$ on A, and so $A \subset \PQ(g_0)$.   Thus, for all $b \in B$, (\QPOStwo) would give
$$g_0(b) = \Phi_A(b) = \supn_A\big[\xbra\cdot{b} - q\big] = \supn_A\big[\xbra\cdot{b} - g_0\big] \le \supn_{\PQ(g_0)}\big[\xbra\cdot{b} - g_0\big].\meqno\EEone$$
Now let $H$ be a nonzero Hilbert space, and consider the SSDB space $\big(H \times H,\xbra\cdot\cdot,\|\cdot\|\big)$, as described at the beginning of this section.   Then\quad $(x,x^*) \in \PQ(g_0) \iff x^* = x$.\quad Let $y \in H\setminus\{0\}$.   Then (\EEone) would give
$$\|y\|^2 = g_0(y,-y) \le \supn_{x \in H}\big[\Xbra{(x,x)}{(y,-y)} - g_0(x,x)\big] = \supn_{x \in H}\big[{-}g_0(x,x)\big] \le 0.$$
Since this is manifestly impossible, there cannot exist a $q$--positive subset $A$ of $B$ such that $\Phi_A = g_0$.   This example is an extension of \cite\HBM, Example 19.20, p.\ 85\endcite.
\defSection \MONsec
\bigbreak
\centerline{\bf \MONsec\quad Monotone multifunctions on reflexive Banach spaces}
\medskip
\noindent
Let $E$ be a nonzero reflexive Banach space.   If $S\colon\ E \toto E^*$ is a multifunction then we use the standard notation
$G(S) := \big\{(x,x^*) \in E \times E^*\colon\ x^* \in Sx\big\}$.   We also write $D(S) := \break\big\{x \in E\colon\ Sx \ne \emptyset\big\} = \pi_1G(S)$.   If $S$ is monotone and $G(S)$ is nonempty, we define the \slant Fitzpatrick functions\endslant, $\phi_S$, associated with $S$ by
$$\phi_S(x,x^*) := \Phi_{G(S)}(x,x^*) = \supn_{(s,s^*) \in G(S)}\big[\bra{x}{s^*} + \bra{s}{x^*} - \bra{s}{s^*}\big].$$
The function $\phi_A$ was introduced by Fitzpatrick in \cite\FITZ, Definition 3.1, p.\ 61\endcite.   The following result was first proved in \cite\FITZ, Corollary 3.9, p.\ 62\endcite\ and  \cite\FITZ, Proposition 4.2, pp.\ 63--64\endcite.
\defTheorem \MAXthm
\medbreak
\noindent
{\bf Theorem  \MAXthm.}\enspace\slant $E$ be a nonzero reflexive Banach space and $S\colon\ E \toto E^*$ be maximally monotone.   Then\quad $\phi_S$ is a BC--function\quand $\PQ\big({\phi_S}^@\big) = \PQ\big(\phi_S\big) = G(S)$.\endslant
\Proof This is immediate from Theorem \PHIMAXthm.\qed 
\medbreak
The following result will be useful in Corollary \SUMcor.
\defLemma \DDOMlem
\medbreak
\noindent
{\bf Lemma \DDOMlem.}\enspace\slant $E$ be a nonzero reflexive Banach space and $S\colon\ E \toto E^*$ be maximally monotone.   Then $D(S) \subset \pi_1\dom\,\phi_S$.\endslant
\Proof From Theorem \MAXthm\ and the finite--valuedness of $q$, $G(S) = \PQ(\phi_S) \subset \dom\,\phi_S$, from which $D(S) = \pi_1G(S) \subset \pi_1\dom\,\phi_S$.\qed
\medbreak
If $S,T\colon\ E \toto F$ then, for all $x \in E$, we use the following standard notation:
$$(S + T)x := \{y + z\colon\ y \in Sx,\ z \in Tx\}.$$
\par
Theorem \SURthm\ is ``Rockafellar's surjectivity theorem'' --- see \cite\SUMS, Proposition 1, p.\ 77\endcite\ for the original proof depending ultimately on Brouwer's fixed--point theorem and an Asplund renorming.   The proof given here is a simplification of that given in Simons--Z\u alinescu, \cite\SZNZ, Theorem 3.1(b), p.\ 8\endcite, and appeared in Simons, \cite\HBM, Theorem 29.6, p.\ 119\endcite.   These references also provide a number of formulae for the exact value of $\min\big\{\|x\|\colon\ x \in E,\ (S + J)x \ni 0\big\}$ in terms of $\phi_{G(S)}$.   Here is one of them:
$$\min\big\{\|x\|\colon\ x \in E,\ (S + J)x \ni 0\big\} = \rthalf \supn_{b \in E \times E^*}\Big[\|b\| - \sqrt{2\phi_S(b) + \|b\|^2}\Big]\vee 0.$$\par
%
\defTheorem \SURthm
\medbreak\noindent
{\bf Theorem \SURthm.}\enspace\slant Let $E$ be a nonzero reflexive Banach space and $S\colon\ E \toto E^*$ be maximally monotone.   Then\quad $(S + J)(E) = E^*$.\endslant
\Proof Let $y^*$ be an arbitrary element of $E^*$.   From Theorem \EETRANSthm(b) \big(with $A := G(S)$\big),\quad $(0,y^*) \in G(S) - G(-J)$.\quad Thus there exist $(s,s^*) \in G(S)$ and $(x,x^*) \in G(J)$ such that\quad $(0,y^*) = (s,s^*) - (x,-x^*)$.\quad But then $x = s$, and so
$$y^* = s^* + x^* \in Sx + Jx = (S + J)x \subset (S + J)(E).\eqno\qed$$
\par
\defRemark \SURrem
\medbreak
\noindent
{\bf Remark \SURrem.}\enspace We note that if $J$ and $J^{-1}$ are single--valued then the converse of Theorem \SURthm\ holds, that is to say\quad $(S + J)(E) = E^* \lr S$ is maximally monotone,\quad however this fails in general.   See Simons, \cite\HBM, Remark 29.7, pp. 120--121\endcite.   It is interesting to observe that, even when the converse of Theorem \SURthm\ fails, the necessary and sufficient condition for maximal monotonicity given in Theorem \EETRANSthm(b) remains true.
\medbreak
We now come to the ``sum theorem''.  More general results can be found in Bo\c t--Csetnek--Wanka, \cite\BCW\endcite\ and Bo\c t--Grad--Wanka, \cite\BGW\endcite.
\defTheorem \SUMthm
\medbreak
\noindent
{\bf Theorem \SUMthm.}\enspace\slant Let $E$ be a nonzero reflexive Banach space, $S,T\colon\ E \toto E^*$ be maximally monotone and $\bigcup_{\lambda > 0}\lambda\big[\pi_1\dom\,\phi_S - \pi_1\dom\,\phi_T\big]$ be a closed linear subspace of $E$.   Then $S + T$ is maximally monotone.\endslant
\Proof From Theorem \MAXthm, $\phi_S$ and $\phi_T$ are BC--functions, $\PQ\big(\phi_S\big) = G(S)$ and $\PQ\big(\phi_T\big) = G(T)$.   Theorem \COMBthm\ with $f := \phi_S$ and $g := \phi_T$ now implies that
$$\big\{(x,s^* + t^*)\colon\ (x,s^*) \in G(S),\ (x,t^*) \in G(T)\big\} = G(S + T)$$
is a maximally monotone subset of $E \times E^*$.   This completes the proof of Theorem \SUMthm.\qed 
\medbreak
The following consequence of Theorem \SUMthm\ is useful in applications:
\defCorollary \SUMcor
\medbreak
\noindent
{\bf Corollary \SUMcor.}\enspace\slant Let $E$ be a nonzero reflexive Banach space, $S,T\colon\ E \toto E^*$ be maximally monotone and $\bigcup_{\lambda > 0}\lambda\big[D(S) - D(T)\big] = E$.   Then $S + T$ is maximally monotone.\endslant
\Proof This is immediate from Lemma \DDOMlem\ and Theorem \SUMthm.\qed
\medbreak
Our next result is ``Rockafellar's sum theorem''.   It follows immediately from Corollary \SUMcor. 
\defCorollary \RTRSUMcor
\medbreak
\noindent
{\bf Corollary \RTRSUMcor.}\enspace\slant Let $E$ be a nonzero reflexive Banach space, $S,T\colon\ E \toto E^*$ be maximally monotone and $D(S) \cap\intr\,D(T) \ne \emptyset$.   Then $S + T$ is maximally monotone.\endslant
\medbreak
A number of other sufficient conditions have been given for the sum of maximally monotone multifunctions on a reflexive Banach space to be maximally monotone.   Many of them are contained in the following `` Sandwiched closed subspace theorem'', which first appeared in Simons--Z\u{a}linescu, \cite\SZNZ, Theorem 5.5, p.\ 13\endcite.
\defTheorem \CQthm
\medbreak
\noindent
{\bf Theorem \CQthm.}\enspace\slant Let $E$ be a nonzero reflexive Banach space, $S,T\colon\ E \toto E^*$ be maximally monotone, and suppose that there exists a closed linear subspace $F$ of  $E$ such that
$$D(S) - D(T) \subset F \subset \bigcupn_{\lambda > 0}\lambda\big[\pi_1\dom\,\phi_S - \pi_1\dom\,\phi_T\big].$$
Then $S + T$ is maximally monotone.   Furthermore, for all $\eps > 0$,
$$D(S) - D(T) \subset \pi_1\dom\,\phi_S - \pi_1\dom\,\phi_T \subset (1 + \eps)\big(D(S) - D(T)\big),$$
(that is to say, the sets\quad $\pi_1\dom\,\phi_S - \pi_1\dom\,\phi_T$ and $D(S) - D(T)$\quad are almost identical) and
$$\bigcupn_{\lambda > 0}\lambda\big[\pi_1\dom\,\phi_S - \pi_1\dom\,\phi_T\big] = \bigcup\nolimits_{\lambda >
0}\lambda\big[D(S) - D(T)\big].$$\endslant\par
\defSection \HAMMsec
\bigbreak
\centerline{\bf \HAMMsec\quad An abstract Hammerstein theorem}
\medskip
Let $E$ be a nonzero reflexive Banach space.   We define the \slant reflection maps\endslant\ $\rho_1,\rho_2$ on $E \times E^*$ by $\rho_1(x,x^*) := (-x,x^*)$ and $\rho_2(x,x^*) := (x,-x^*)$.   Our first result is the ``$\rho_1$--transversality theorem''.
\defTheorem \RHOTRANSthm
\medbreak
\noindent
{\bf Theorem \RHOTRANSthm.}\enspace\slant Let $E$ be a nonzero reflexive Banach space, $f,g \in \PCLSC(E \times E^*)$ and $f,g$ be BC--functions.   Then:
$$\dom\,f + \rho_1\dom\,g  = E \times E^* \iff \PQ(f) + \rho_1\PQ(g) = E \times E^*$$
and
$$\dom\,f - \rho_1\dom\,g  = E \times E^* \iff \PQ(f) - \rho_1\PQ(g) = E \times E^*.$$\endslant\par
\Proof This follows from Corollary \ABPOSNEGcor\ with $\rho := \rho_2$, and the fact that $\rho_1 = -\rho_2^{-1}$.\qed 
\medbreak
We give the next two results for completeness, though more general results are known.  See, for instance, Z\u alinescu, \cite\ZFITZ, Theorem 3 and Corollary 4\endcite. 
\defTheorem \RHOSURthm
\medbreak
\noindent
{\bf Theorem \RHOSURthm.}\enspace\slant Let $E$ be a nonzero reflexive Banach space, $S,T\colon\ E \toto E^*$ be maximally monotone and \quad $\dom\,\phi_S + \rho_1\dom\,\phi_T  = E \times E^*$.\quad  Then\quad $(S + T)(E) = E^*$.\endslant
\Proof Let $y^*$ be an arbitrary element of $E^*$.  From Theorem \RHOTRANSthm\ (with $f := \phi_S$ and $g := \phi_T$) and Theorem \MAXthm,\quad  $(0,y^*) \in G(S) + \rho_1G(T)$.\quad Thus there exist $(s,s^*) \in G(S)$ and $(x,x^*) \in G(T)$ such that\quad $(0,y^*) = (s,s^*) + (-x,x^*)$.\quad But then $x = s$, and so
$$y^* = s^* + x^* \in Sx + Tx = (S + T)x \subset (S + T)(E).\eqno\qed$$
\par
Since\quad $G(J) = \PQ(g_0)$,\quad it follows that\quad $\phi_J = \Phi_{G(J)} = \Phi_{\PQ(g_0)}$.\quad Thus Remark \Grem\ implies that\quad $\dom\,\phi_J = E \times E^*$.\quad   Consequently, the following result generalizes Theorem \SURthm.   Compare with \Cite\MLSUR, Theorem 2\big((1)$\lr$(2)\big)\endcite. 
\defCorollary \RHOSURcor
\medbreak
\noindent
{\bf Corollary \RHOSURcor.}\enspace\slant Let $E$ be a nonzero reflexive Banach space, $S,T\colon\ E \toto E^*$ be maximally monotone and \quad $\dom\,\phi_T  = E \times E^*$.\quad  Then\quad $(S + T)(E) = E^*$.\endslant
\Proof This is immediate from Theorem \RHOSURthm.\qed
\defTheorem \RRABthm
\medbreak
\noindent
{\bf Theorem \RRABthm.}\enspace\slant Let $E$ be a nonzero reflexive Banach space, $f,g \in \PCLSC(E \times E^*)$ be BC--functions, and $w^* \in E^*$ be such that $E \times \{w^*\} \subset \dom\,f$ and $\pi_2\,\dom\,g = E^*$.   Then:
\smallbreak
\noindent
{\rm(a)}\enspace $\PQ(f) + \rho_1\,\PQ(g) = E \times E^*$ and $\PQ(f) - \rho_1\,\PQ(g) = E \times E^*$.
\smallbreak
\noindent
{\rm(b)}\enspace If $x \in E$ then there exist $(y,y^*) \in \PQ(f)$ and $(z,y^*) \in \PQ(g)$ such that $y + z = x$.
\smallbreak
\noindent
{\rm(c)}\enspace If $x^* \in E^*$ then there exist $(y,y^*) \in \PQ(f)$ and $(y,z^*) \in \PQ(g)$ such that $y^* + z^* = x^*$.
\endslant
\Proof(a)\enspace Let $(x,x^*)$ be an arbitrary element of $E \times E^*$.   Since\quad $\pi_2\,\dom\,g = E^*$,\quad there exist\quad $y,z \in E$\quad such that\quad $(y,x^* - w^*) \in \dom\,g$\quad and\quad $(z,w^* - x^*) \in \dom\,g$.\quad   But\quad $(x + y,w^*) \in \dom\,f$\quad and\quad $(x - z,w^*) \in \dom\,f$,\quad hence
$$(x,x^*) = (x + y,w^*) + \rho_1(y,x^* - w^*) \in \dom\,f + \rho_1\,\dom\,g$$
and
$$(x,x^*) = (x - z,w^*) - \rho_1(z,w^* - x^*) \in \dom\,f - \rho_1\,\dom\,g.$$
Thus we have proved that\quad$\dom\,f + \rho_1\,\dom\,g = E \times E^*$\quand $\dom\,f - \rho_1\,\dom\,g = E \times E^*$, and (a) follows from Theorem \RHOTRANSthm.
\smallbreak
(b)\enspace It follows from (a) that there exist\quad $(y,y^*) \in \PQ(f)$\quad and\quad $(z,z^*) \in \PQ(g)$\quad such that\quad $(y,y^*) + (z,-z^*) = (x,0)$.\quad   But then\quad $z^* = y^*$\quad and\quad $y + z = x$. 
\smallbreak
(c)\enspace It follows from (a) that there exist\quad $(y,y^*) \in \PQ(f)$\quad and\quad $(z,z^*) \in \PQ(g)$\quad such that\quad $(y,y^*) - (z,-z^*) = (0,x^*)$.\quad   But then\quad $z = y$\quad and\quad $y^* + z^* = x^*$.\qed
\medbreak
We now reverse the direction of $T$.
\defTheorem \STINVthm
\medbreak
\noindent
{\bf Theorem \STINVthm.}\enspace\slant Let $E$ be a nonzero reflexive Banach space and $S\colon\ E \toto E^*$ and $T\colon\ E^* \toto E$ be maximally monotone.   Suppose that $\pi_1\,\dom\,\phi_T = E^*$ and there exists $w^* \in E^*$ such that $E \times \{w^*\} \subset \dom\,\phi_S$.   Then:
\smallbreak
\noindent{\rm(a)}\enspace If $I_E$ is the identity map on $E$, $(I_E + TS)(E) = E$.
\smallbreak
\noindent{\rm(b)}\enspace If $I_{E^*}$ is the identity map on $E^*$, $(I_{E^*} + ST)(E^*) = E^*$.\endslant
\Proof Let $f := \phi_S$ and $g := \phi_{T^{-1}}$, so that $\pi_2\,\dom\,g = E^*$.\smallbreak
(a)\enspace Let $x$ be an arbitrary element of $E$.   From Theorem \RRABthm(b) and Theorem \MAXthm, there exist \quad $(y,y^*) \in G(S)$\quad and \quad $(y^*,z) \in G(T)$\quad such that \quad $y + z = x$. \quad   Thus\quad $z \in Ty^* \subset TSy$\quad and\quad $x = y + z \in (I_E + TS)(y) \subset (I_E + TS)(E)$.
\smallbreak
(b)\enspace Let $x^*$ be an arbitrary element of $E^*$.   From Theorem \RRABthm(c) and Theorem \MAXthm, there exist\quad $(y,y^*) \in G(S)$\quad and \quad $(z^*,y) \in G(T)$\quad such that \quad $y^* + z^* = x^*$.   Thus\quad $y^* \in Sy \subset STz^*$\quad and \quad $x^* = z^* + y^* \in (I_{E^*} + ST)(z^*) \subset (I_{E^*} + ST)(E^*)$.\qed
\medbreak
The next result is a considerable generalization of \cite\ZEIDLER, Theorem 32.O, p.\ 909\endcite\ which, in turn, was applied to Hammerstein integral equations.   See \cite\HBM, Remark 30.5, p.\ 124\endcite\ for a more complete discussion.
\defTheorem \HAMMERthm
\medbreak
\noindent
{\bf Theorem \HAMMERthm.}\enspace\slant Let $E$ be a nonzero reflexive Banach space and $S\colon\ E \toto E^*$ and $T\colon\ E^* \toto E$ be maximal monotone.   Suppose that \quad {\bf either}\quad $\pi_1\,\dom\,\phi_T = E^*$\quad and there exists\break $w^* \in E^*$\quad such that \quad $E \times \{w^*\} \subset \dom\,\phi_S$\quad {\bf or} \quad $\pi_1\,\dom\,\phi_S = E$\quad and there exists \quad $w \in E$\quad such that \quad $E \times \{w\} \subset \dom\,\phi_T$. \quad Then $(I_E + TS)(E) = E$.
\endslant
\Proof The first case has already been established in Theorem \STINVthm(a), while the\break second case follows from Theorem \STINVthm(b), with $E$ replaced by $E^*$, and the roles of $S$ and $T$ interchanged.\qed
\defRemark \FIFTYrem
\medbreak
\noindent
{\bf Remark \FIFTYrem.}\enspace In the recent paper \cite\FIFTY\endcite, Jonathan Borwein has written a survey of the history of monotonicity (even in the nonreflexive case) over the past fifty years.
\bigbreak
\centerline{\bf References}
\medskip
\nmbr\AB
\item{[\AB]} H. Attouch and H. Brezis, {\sl  Duality for the sum of
convex funtions in general Banach spaces.}, Aspects of Mathematics and its Applications, J. A. Barroso, ed., Elsevier Science Publishers (1986), 125--133.
\nmbr\FIFTY
\item{[\FIFTY]} J. M. Borwein, \slant Fifty years of maximal monotonicity\endslant, $<$http://www.carma.newcastle. edu.au/$\sim$jb616/fifty.pdf$>$. 
\nmbr\BCW
\item{[\BCW]}R. I. Bo\c t, E. R. Csetnek and G. Wanka, \slant A new condition for maximal monotonicity via representative functions\endslant, Nonlinear Anal. {\bf 67} (2007), 2390--2402.
\nmbr\BGW
\item{[\BGW]}R. I. Bo\c t, S--M Grad and G. Wanka, \slant Maximal monotonicity for the precomposition with a linear operator\endslant.   SIAM J. Optim. {\bf 17} (2006), 1239--1252. 
\nmbr\BS
\item{[\BS]} R. S. Burachik and B. F. Svaiter, \slant Maximal monotonicity, conjugation and the duality product\endslant,  Proc. Amer. Math. Soc.  {\bf 131} (2003), 2379--2383.
\nmbr\FITZ
\item{[\FITZ]} S. Fitzpatrick, \slant Representing monotone operators
by convex functions\endslant,   Workshop/\break Miniconference on
Functional Analysis and Optimization (Canberra, 1988),  59--65, Proc.
Centre Math. Anal. Austral. Nat. Univ., {\bf 20}, Austral. Nat. Univ.,
Canberra, 1988.
\nmbr\MLSUR
\item{[\MLSUR]} J-E. Mart\'\i nez-Legaz, \slant Some generalizations of Rockafellar's surjectivity theorem\endslant,\break Pacific J. of Optimization {\bf 4} (2008), 527--535.
\nmbr\MLT
\item{[\MLT]} J.--E. Mart\'\i nez-Legaz and M. Th\'era, \slant
A convex representation of maximal monotone operators\endslant,\
J. Nonlinear Convex Anal. {\bf 2} (2001), 243--247.
\nmbr\PENOT
\item{[\PENOT]}J.--P. Penot, \slant The relevance of convex analysis for the study of monotonicity\endslant,  Nonlinear Anal. {\bf 58}  (2004), 855--871.
\nmbr\PZA
\item{[\PZA]}J.--P. Penot and C. Z\u{a}linescu, \slant Some problems about the representation of monotone operators by convex functions\endslant,  ANZIAM J. {\bf 47}  (2005), 1--20.
\nmbr\RTRFENCHEL
\item{[\RTRFENCHEL]} R. T. Rockafellar, \slant Extension of Fenchel's
duality theorem for convex functions\endslant, Duke Math. J. {\bf33}
(1966),  81--89.
\nmbr\SUMS
\item{[\SUMS]}-----, {\sl On the Maximality of Sums
of Nonlinear Monotone Operators}, Trans. Amer. Math. Soc. {\bf
149}(1970),75-88.
\nmbr\MANDM
\item{[\MANDM]}S. Simons, \slant Minimax and monotonicity\endslant, 
Lecture Notes in Mathematics {\bf 1693} (1998),\break Springer--Verlag.
\nmbr\HBM
\item{[\HBM]}-----, \slant From Hahn--Banach to monotonicity\endslant, 
Lecture Notes in Mathematics, {\bf 1693},\break second edition, (2008), Springer--Verlag.
\nmbr\SSDMON
\item{[\SSDMON]}-----, \slant Banach SSD spaces and classes of monotone sets\endslant, http://arxiv/org/abs/\break 0908.0383v2, posted August 26, 2009.
\nmbr\SZNZ
\item{[\SZNZ]}S. Simons and C. Z\u{a}linescu, \slant Fenchel duality,
Fitzpatrick functions and maximal monotonicity\endslant, J. of
Nonlinear and  Convex Anal., {\bf 6} (2005), 1--22.
\nmbr\VZ
\item{[\VZ]}M. D. Voisei and C. Z\u{a}linescu, \slant Strongly--representable operators\endslant, J. of Convex Anal., {\bf 16(3)}, 2009. 
\nmbr\ZBOOK
\item{[\ZBOOK]} C. Z\u{a}linescu, \slant Convex analysis in
general vector spaces\endslant, (2002), World Scientific.
\nmbr\ZFITZ
\item{[\ZFITZ]} -----, \slant A new convexity property for monotone operators\endslant, J. Convex Anal. {\bf 13} (2006), 883--887.
\nmbr\ZEIDLER
\item{[\ZEIDLER]} E. Zeidler, \slant Nonlinear functional analysis and its applications, II/B, Nonlinear monotone operators\endslant, Springer-Verlag, New York, 1990.
\Signoff
\bigskip
Department of Mathematics\par
University of California\par
Santa Barbara\par
CA 93106-3080\par
U. S. A.\par
email:  simons@math.ucsb.edu
\bye